\documentclass{gtpart}     
\usepackage{amssymb}
\usepackage{xypic}
\usepackage{graphicx}
\usepackage{amsmath}

\title{The $R(S^1)$--graded equivariant homotopy of THH($\mathbb{F}_p$)}

\author{Teena Gerhardt}
\givenname{}
\surname{}
\address{Indiana University Bloomington}
\email{tgerhard@indiana.edu}
\urladdr{http://mypage.iu.edu/~tgerhard}

\keyword{}
\subject{primary}{msc2000}{19D55}
\subject{secondary}{msc2000}{55P91}

\arxivreference{0710.2938}
\arxivpassword{xue8h}

\volumenumber{}
\issuenumber{}
\publicationyear{}
\papernumber{}
\startpage{}
\endpage{}
\doi{}
\MR{}
\Zbl{}
\received{}
\revised{}
\accepted{}
\published{}
\publishedonline{}
\proposed{}
\seconded{}
\corresponding{}
\editor{}
\version{}

\newtheorem{thm}{Theorem}[section]    
\newtheorem{lem}[thm]{Lemma}          
\theoremstyle{definition}

\newtheorem{proposition}[thm]{Proposition}
\newtheorem{corollary}[thm]{Corollary}

\newcommand{\TR}{\mathrm{TR}}
\newcommand{\TC}{\mathrm{TC}}

\begin{document}

\begin{abstract}

The main result of this paper is the computation
of TR$^n_{\alpha}(\mathbb{F}_p;p)$ for $\alpha \in R(S^1)$.
These $R(S^1)$--graded TR--groups are the equivariant homotopy
groups naturally associated to the $S^1$--spectrum
THH$(\mathbb{F}_p)$, the topological Hochschild $S^1$--spectrum.
This computation, which extends a partial result of Hesselholt and
Madsen, provides the first example of the
$R(S^1)$--graded TR--groups of a ring. These groups arise in
algebraic $K$--theory computations, and are particularly important
to the understanding of the algebraic $K$--theory of non-regular
schemes.

\end{abstract}

\maketitle


\section{Introduction}
Algebraic $K$--theory groups are generally very difficult to compute, but B\"{o}kstedt, Hsiang,
and
Madsen \cite{Trace} developed a homotopy-theoretic approach to algebraic $K$--theory computations  that has been quite fruitful. For every ring $A$, they defined a cyclotomic trace map
\begin{displaymath}
\mathrm{trc}: K_q(A) \rightarrow \mathrm{TC}_q(A)
\end{displaymath}
relating algebraic $K$--theory to topological cyclic homology. This map is highly non-trivial, and thus one can often understand algebraic
$K$--theory by understanding topological cyclic homology (see,
for instance, McCarthy \cite{McCarthy}).

As an approach to computing topological cyclic homology, Hesselholt and Madsen studied TR--theory ~\cite{KWitt, local}. For a ring $A$ and a
fixed prime $p$,
\begin{displaymath}
\mathrm{TR}^n_q(A;p) :=
\pi_q(\mathrm{T}(A)^{C_{p^{n-1}}}) = [S^q \wedge
{S^1/C_{p^{n-1}}}_+, \mathrm{T}(A)]_{S^1}
\end{displaymath}
where T$(A)$ denotes the topological Hochschild $S^1$--spectrum of
$A$ and T$(A)^{C_{p^{n-1}}}$ denotes the $C_{p^{n-1}}$ fixed
point spectrum. These TR--groups come equipped with
several operators and relations which provide a rigid algebraic structure, making computations possible. Topological cyclic homology is defined by a
homotopy limit construction
that involves these operators. Thus understanding the TR--groups of
a ring helps us to understand its topological cyclic homology, and
hence its algebraic $K$--theory.

The topological Hochschild $S^1$--spectrum also has naturally associated equivariant homotopy groups which give a TR--theory graded by the complex representation ring of the circle, $R(S^1)$. These groups arise naturally from the computational method outlined above. Elements in the representation ring are given by formal differences of isomorphism classes of representations. For every $\alpha \in R(S^1)$ we choose representatives $\beta$ and $\gamma$ such that $\alpha = [\beta] -[\gamma]$. Let $S^{\beta}$ denote the one--point compactification of the representation $\beta$. Then the $R(S^1)$--graded TR--groups are defined as
\begin{displaymath}
\TR^n_{\alpha}(A;p)  = [S^{\beta} \wedge {S^1/C_{p^{n-1}}}_+, S^{\gamma} \wedge T(A) ]_{S^1}.
\end{displaymath}

These $R(S^1)$--graded groups first arose in computations of the algebraic $K$--theory of non-regular schemes.  For instance, for an $\mathbb{F}_p$--algebra $A$, Hesselholt and Madsen ~\cite{Poly} expressed  the $K$--theory of $A[x]/(x^e)$ in terms of the $R(S^1)$--graded $\TR$--groups of A. While $R(S^1)$--graded TR--groups have been used in computations, there are no fully computed examples for all $\alpha \in R(S^1)$. The first computation to be done is that of $\TR^n_{\alpha}(\mathbb{F}_p; p)$, which is the main result of this paper.

Before stating the result of this computation, we introduce some notation. Let
$$
\rho_p: S^1 \rightarrow S^1/C_p
$$
be the isomorphism given by the $p$th root. Then we define a prime operation as follows: for $\alpha \in R(S^1)$,
$$
\alpha' = \rho_p^*(\alpha^{C_p}).
$$
Recall the ring structure of $R(S^1)$:
$$
R(S^1) \cong \mathbb{Z}[t, t^{-1}].
$$
If $p$ divides $i$ then $(t^i)' = t^{i/p}$. If $p$ does not divide $i$ then $(t^i)' = 0$. This prime operation is additive. We use the notation $\alpha^{(k)}$ to denote the prime operation applied $k$ times to $\alpha$. We denote the dimension of $\alpha$ as a complex vector space by $|\alpha|$. Note that every finite sequence of integers can be realized as the dimensions $|\alpha|, |\alpha'|, \ldots, |\alpha^{(n)}|$ for some virtual representation $\alpha \in R(S^1)$.

Our main result expresses $\TR^{n}_{\alpha}(\mathbb{F}_p;p)$ as an iterated pullback diagram that can be computed explicitly.  The argument is inductive, computing the group $\TR^{n}_{ \alpha}(\mathbb{F}_p;p)$ from the group $\TR^{n-1}_{\alpha'}(\mathbb{F}_p;p)$ using the  fundamental long exact sequence of TR--theory:
$$
\xymatrix@C=9pt{
\cdots \ar[r] &  \pi_{\alpha}T(\mathbb{F}_p)_{hC_{p^{n-1}}}\ar[r]  & \TR^n_{\alpha}(\mathbb{F}_p;p)  \ar[r]  & \TR^{n-1} _{\alpha '} (\mathbb{F}_p;p)  \ar[r]  & \cdots &\\}
$$
In order to determine this extension we use information about a map
$$
\hat{\Gamma}_{n-1}: \TR^{n-1}_{\alpha'}(\mathbb{F}_p;p) \rightarrow
\pi_{\alpha}\hat{\mathbb{H}}(C_{p_{n-1}}, T(\mathbb{F}_p))
$$
fitting into a diagram of long exact sequences
$$
\xymatrix@C=9pt{
\cdots \ar[r]&  \pi_{\alpha}T(\mathbb{F}_p)_{hC_{p^{n-1}}}  \ar[r] \ar@{=}[d] & \TR^n_{\alpha}(\mathbb{F}_p;p)  \ar[r] \ar[d]^{\Gamma_{n-1}} & \TR^{n-1} _{\alpha '} (\mathbb{F}_p;p)  \ar[r] \ar[d]^{\hat{\Gamma}_{n-1}} & \cdots \\
\cdots \ar[r] &  \pi_{\alpha}T(\mathbb{F}_p)_{hC_{p^{n-1}}} \ar[r] &
 \pi_{\alpha}T(\mathbb{F}_p)^{hC_{p^{n-1}}} \ar[r] &
 \pi_{\alpha}\hat{\mathbb{H}}(C_{p_{n-1}}, T(\mathbb{F}_p)) \ar[r]& \cdots
\\}
$$
Here $\hat{\mathbb{H}}(C_{p_{n-1}}, T(\mathbb{F}_p))$ denotes the Tate spectrum. Section \ref{Spectral} contains a more detailed discussion of this diagram of long exact sequences.  We have spectral sequences which allow us to compute the groups in the bottom row of this diagram. We use these diagrams of long exact sequences to identify $\TR^n_{\alpha}(\mathbb{F}_p;p)$ as an iterated pullback.

We now state our main result

\begin{thm}\label{main}
The group $\TR^n_{\alpha}(\mathbb{F}_p;p)$ is given by the following iterated pullback. Let $r_j = \min(j,  |\alpha^{(n-j)}| + 1)$.
$$
\TR^1_{\alpha^{(n-1)}}(\mathbb{F}_p;p) \cong \left\{ \begin{array}{cc}
 \mathbb{Z}/p& \textup{ if }r_1 = 1  \\
 0 & \textup{  else  }
 \end{array}\right.
$$
The map $\hat{\Gamma}_1: \TR^1_{\alpha^{(n-1)}}(\mathbb{F}_p;p) \rightarrow \mathbb{Z}/p$ is injective. Letting $j$ range from 2 to n, if $r_j \geq 1$, $\TR^j_{ \alpha^{(n-j)}}(\mathbb{F}_p;p)$ is given by the following pullback diagram, which also uniquely determines the map $\hat{\Gamma}_j$ up to units:
$$
\xymatrix{
\TR^j_{\alpha^{(n-j)}}(\mathbb{F}_p;p) \ar[r]^{R_j} \ar[d]^{\hat{\Gamma}_j} &\TR^{j-1}_{\alpha^{(n-j+1)}}(\mathbb{F}_p;p) \ar[d]^{\hat{\Gamma}_{j-1}} \\
\mathbb{Z}/p^j \ar[r]^{R^h_j} & \mathbb{Z}/p^{j-1} \\}
$$
where $R^h_j(1) \doteq p^{r_j - 1}$. If $r_j <1$, $\TR^j_{ \alpha^{(n-j)}}(\mathbb{F}_p;p) \cong \TR^{j-1}_{ \alpha^{(n-j+1)}}(\mathbb{F}_p;p) $ and $\hat{\Gamma}_j$ is the composite $\phi_j\hat{\Gamma}_{j-1}$ where $\phi_j:\mathbb{Z}/p^{j-1} \rightarrow \mathbb{Z}/p^{j}$, $\phi_j(1) \doteq p^{1-r_j}$.
\end{thm}

The following theorem can be proven by computing these pullbacks
\begin{thm}\label{main2}
There is an isomorphism of abelian groups
$$
\xymatrix{ \mathbb{Z}/p^{l_{1,n}} \oplus \mathbb{Z}/p^{l_{2,n}} \oplus \cdots \oplus \mathbb{Z}/p^{l_{n,n}} \ar[r]^-{\iota_n}_-{\sim} & \TR^n_{\alpha}(\mathbb{F}_p;p) }
$$
where the $l_{i,n}$ are the explicitly defined recursive functions of  $|\alpha|, |\alpha'|, \ldots, |\alpha^{(n-1)}|$ given in Section \ref{Proofs}.
\end{thm}

We make a few observations about the results of these computations. If the integers $|\alpha|, |\alpha'|, \ldots |\alpha^{(n-1)}|$ are all negative, then $\TR^n_{\alpha}(\mathbb{F}_p;p) = 0.$ If the integers $|\alpha|, |\alpha'|, \ldots |\alpha^{(n-1)}|$ are all nonnegative, then $\TR^n_{\alpha}(\mathbb{F}_p;p) \cong \mathbb{Z}/p^n$. Thus the TR--groups are easily understood for certain representations $\alpha$. For any $\alpha$, each of the values $l_{i,n}$ given in the computation of $\TR^n_{\alpha}(\mathbb{F}_p;p)$ satisfies $l_{i,n} \leq n$ since $TR^n_{* + \alpha}(\mathbb{F}_p;p)$ is a $TR^n_*(\mathbb{F}_p;p)$--module, where $*$ denotes an even integer grading. However, the sum $l(n) = l_{1, n} + l_{2,n} + \cdots + l_{n,n}$ can be greater than $n$.

These computations extend a result of Hesselholt and Madsen \cite[Proposition 9.1]{KWitt} for representations of the form $\alpha = q - \gamma$, $q \in \mathbb{Z}$, $\gamma$ a complex representation of $S^1$. If the virtual representation $\alpha$ is of the form $q + \lambda$ or $q-\lambda$ for an actual complex $S^1$--representation $\lambda$, then the groups $\TR^n_{\alpha}(\mathbb{F}_p;p)$ are cyclic. The computations of the TR--groups in these special cases is presented in Section \ref{Proofs}.

In algebraic $K$--theory computations, it is the $RO(S^1)$--graded equivariant homotopy groups, rather than the $R(S^1)$--graded groups, that arise most naturally. The even dimensional real representations correspond to the complex representations, the computations for which appear in this paper. The groups graded by odd--dimensional real representations must be computed slightly differently. These computations will appear in joint work with Vigleik Angeltveit.

\subsection{Organization}
In Section \ref{TRSection} we discuss a homotopy-theoretic approach to algebraic $K$--theory computations. In particular, we recall the definition of TR--theory, and the operators and relations on it.  We then define the $R(S^1)$--graded TR--groups of a ring $A$. In Section \ref{Spectral} we study the fundamental long exact sequence of TR--theory, and its relation to the norm--restriction sequence. Section \ref{compute} focuses on applying this relationship to execute the first half of the induction step needed to prove Theorems  \ref{main} and \ref{main2}.  In Section \ref{HatGamma} we finish the induction step by studying the map $\hat{\Gamma}$.  In Section \ref{Proofs} we put our earlier results together to prove Theorems \ref{main} and \ref{main2} and several corollaries.

\subsection{Notation and conventions}

Throughout this paper, $A$ denotes a commutative ring, and $p$ a fixed prime. The symbol $\alpha$ denotes a virtual representation in $R(S^1)$. We use the convention $a \doteq b$ for $a=\mu b$, where $\mu$ is a unit.

\subsection{Acknowledgments}

I would like to thank Lars Hesselholt for all of his support and guidance on this project. I am also grateful to Vigleik Angeltveit for many helpful conversations.

\section{Introduction to TR--groups}\label{TRSection}

\subsection{Overview of the approach}

Let $p$ be a fixed prime and $A$ a commutative ring. In this section we outline a homotopy-theoretic approach to algebraic $K$--theory computations. In particular, we recall the definitions and properties of TR--groups, studied by Hesselholt and Madsen
(see ~\cite{KWitt, local}).

To every ring $A$, there is an associated genuine $S^1$--spectrum, $T(A)$, whose zero'th space is naturally equivalent to the topological Hochschild space THH($A$) defined by B\"{o}kstedt \cite{Bok}.
This topological Hochschild spectrum is cyclotomic \cite{KWitt}, which means in particular that it comes equipped with an equivalence of $S^1$--spectra
$$
r: \rho_p^*(\tilde{E} \wedge T(A))^{C_p}
\rightarrow T(A).
$$
Here $E$ denotes a free contractible $S^1$--space, and $\tilde{E}$ is the cofiber of the map $E_{+} \rightarrow S^0$ given by projection onto the non-basepoint. This equivalence will be essential to the construction of the long exact sequence of TR--groups in Section \ref{Spectral}.

The $S^1$--spectrum structure of $T(A)$ allows us to take
$C_{p^{n}}$--fixed point spectra of $T(A)$ for $C_{p^{n}} \subset
S^1$ the cyclic group of order $p^n$. Then inclusion of fixed
points induces a map
$$
F: T(A)^{C_{p^n}} \rightarrow T(A)^{C_{p^{n-1}}}
$$
called the Frobenius. We also get a restriction map
$$
R: T(A)^{C_{p^{n}}} \rightarrow T(A)^{C_{p^{n-1}}}
$$
which we define in Section \ref{Spectral}. Using these operators, we can define
the topological cyclic homology spectrum of $A$. Let
$$
\TC(A;p) = [\underset{\longleftarrow}{\mathrm{holim}} T(A)^{C_{p^n}}]^{hF}
$$
where the homotopy limit is take across the maps $R$ above, and $X^{hF}$ denotes the $F$--homotopy fixed points of $X$, i.e. the homotopy fiber of id $- F$.

It is possible to repeat the above, working with all natural numbers rather than just powers of one prime. The resulting spectrum TC($A$), however, does not carry much more information than the product of the TC($A; p)$ spectra. In particular, after profinite completion
$$
\TC(A)^{\wedge} \simeq \prod \TC(A;p)^{\wedge}
$$

For every ring A, we have a cyclotomic trace map \cite{Trace} from the algebraic
$K$--theory of $A$ to the topological cyclic homology spectrum of
$A$
\begin{displaymath}
\mathrm{trc}: K(A) \rightarrow \mathrm{TC}(A).
\end{displaymath}
We can often understand algebraic $K$--theory by understanding
topological cyclic homology and this cyclotomic trace map. For example, for $I \subset A$ a nilpotent ideal, McCarthy's theorem \cite{McCarthy} says that after profinite completion, there is an equivalence
$$
\mathrm{trc }: K(A, I) \stackrel{\sim}{\longrightarrow} \mathrm{TC}(A, I),
$$
where $K(A, I)$ and TC($A, I)$ denote relative $K$--theory and relative TC respectively.

\subsection{$\mathbb{Z}$--graded TR--groups}

As an approach to understanding topological cyclic homology,
Hesselholt and Madsen studied TR--groups. For a ring $A$ and a fixed prime
$p$, let
$$
\TR^n(A;p) = T(A)^{C_{p^{n-1}}}
$$
and
\begin{displaymath}
\TR^n_q(A;p) = \pi_q(T(A)^{C_{p^{n-1}}}).
\end{displaymath}
These TR--groups come equipped with several operators. Again,
inclusion of fixed points induces the Frobenius map
\begin{displaymath}
F: \TR^n_q(A;p) \rightarrow \TR^{n-1}_q(A;p).
\end{displaymath}
This map has an associated transfer, the Verschiebung
\begin{displaymath}
V: \TR^{n-1}_q(A;p) \rightarrow \TR^n_q(A;p).
\end{displaymath}
There is also a derivation $d: \TR^n_q(A;p) \rightarrow
\TR^n_{q+1}(A;p)$ induced from the circle action on $T(A)$.
Lastly,
we have a restriction map
\begin{displaymath}
R: \TR^n_q(A;p) \rightarrow \TR^{n-1}_q(A;p)
\end{displaymath}
which we define in Section \ref{Spectral}. These maps satisfy the relations
$FV=p, FdV = d,$ and $VF = V(1)$, where $1$ is the multiplicative unit $[1]_n$ with Witt coordinates $(1,0, \ldots 0) \in W_n(A) \cong \TR_0^n(A;p)$. For a more detailed introduction to TR--theory, see for instance \cite{KWitt, local}.

We can define TC$^n(A;p)$ as the homotopy equalizer of the
maps \begin{displaymath} R, F: \TR^n(A;p) \rightarrow
\TR^{n-1}(A;p).
\end{displaymath}
Then TC$(A;p)$ is the homotopy limit of the spectra TC$^n(A;p)$.
Thus understanding the TR--groups of a ring helps us to
understand its topological cyclic homology, and hence its
algebraic $K$--theory.

Computations are possible on the level of TR because the
operators on TR--groups and the relations between them give
TR--groups a
rigid algebraic structure. Indeed they have the structure of a
Witt
complex over $A$  \cite{deRham}.

\subsection{$R(S^1)$--graded TR--groups}

Recall that TR--groups are defined as
\begin{displaymath}
\TR^n_q(A;p) = \pi_q(T(A)^{C_{p^{n-1}}})
\end{displaymath}
or equivalently
\begin{displaymath}
\TR^n_q(A;p) \cong [S^q \wedge {S^1/C_{p^{n-1}}}_+,  T(A)]_{S^1}.
\end{displaymath}

To any $S^1$--spectrum $T$ there are naturally associated equivariant homotopy groups graded by the representation ring of the circle, $R(S^1)$. Elements in the representation ring are given by formal differences of isomorphism classes of representations. For every $\alpha \in R(S^1)$ we choose representatives $\beta$ and $\gamma$ such that $\alpha = [\beta] -[\gamma]$. Then for an $S^1$--spectrum $T$, the group
$$
\pi_{\alpha}(T) = [S^{\beta}, S^{\gamma} \wedge T]_{S^1},
$$
the abelian group of maps in the $S^1$--stable category. Further, we get an $R(S^1)$--graded Mackey functor $\underline{\pi}_*$ which for each subgroup $C \subset S^1$ assigns
$$
\underline{\pi}_{\alpha}(T)(C) = [S^{\beta} \wedge S^1/C_+, S^{\gamma} \wedge T]_{S^1}
$$
Lewis and Mandell ~\cite{lewis} have proven that this definition gives a well--behaved theory of $R(S^1)$--graded homotopy groups. Although the theorem in the cited reference is stated only for finite groups, the same method of proof works for compact lie groups.

We apply this $R(S^1)$--graded Mackey functor to the topological Hochschild $S^1$--spectrum of $A$ to get $R(S^1)$--graded TR--groups,
\begin{displaymath}
\TR^n_{\alpha}(A;p) = [S^{\beta} \wedge
{S^1/C_{p^{n-1}}}_+, S^{\gamma} \wedge T(A) ]_{S^1}.
\end{displaymath}
These $R(S^1)$--graded TR--groups arise naturally in computations of algebraic $K$--theory.
Given a ring $A$ and a pointed monoid $\Pi$, let $A(\Pi)$ denote the pointed monoid algebra. To study the algebraic $K$--theory of $A(\Pi)$, we start with the topological Hochschild $S^1$--spectrum of $A(\Pi)$. Formal properties of the construction give a well known equivalence
of $S^1$--spectra \cite{Trunc}
$$
T(A) \wedge N^{cy}(\Pi)  \stackrel{\sim}{\rightarrow} T(A(\Pi)).
$$
where $N^{cy}(\Pi)$ denotes the cyclic bar construction on $\Pi$.
For example, to understand the topological Hochschild homology of $A[x]/(x^e)$,
$A$ an $\mathbb{F}_p$ algebra, Hesselholt and Madsen \cite{Poly}
studied
$$
T(A) \wedge N^{cy}(\Pi_e) \simeq T(A(\Pi_e)) = T(A[x]/(x^e)).
$$
Here $\Pi_e$ denotes the pointed multiplicative monoid $\{0, 1, x,
\ldots x^{e-1}\}$ where $x^e=0$. Using the above equivalence, and the definition of TR--groups, we obtain an isomorphism
$$
\textup{TR}^n_q(A[x]/(x^e); p) \cong [S^q \wedge {S^1/C_{p^{n-1}}}_+, T(A) \wedge N^{cy}(\Pi_e)]_{S^1}.
$$
Hesselholt and Madsen demonstrated how the $S^1$--equivariant homotopy type of $N^{cy}(\Pi_e)$ is built from representation spheres. Replacing  $N^{cy}(\Pi_e)$ by these representation spheres in the isomorphism above, we end up in the $RO(S^1)$--graded TR groups of $A$. Further, for virtual representations $\alpha$ of the form arising in this example, the groups $\TR^n_{\alpha}(\mathbb{F}_p;p)$ are zero if the dimension of $\alpha$ as a real vector space is odd. The representations with even real dimension correspond exactly to the complex representations of $\alpha$. So this expresses the algebraic $K$--theory of $A[x]/(x^e)$ in terms of the $R(S^1)$--graded TR--groups of the ring $A$.

The general principle here is that the $R(S^1)$ graded TR--groups of a ring $A$ arise in the computation of the ordinary TR--groups of the more complicated ring $A(\Pi)$.  By considering different pointed monoids $\Pi$, one can use the
$R(S^1)$--graded TR--groups of $A$ to understand the algebraic $K$--theory of a
variety of rings $A(\Pi)$. Using the pointed monoid $\Pi = \{0, 1,
x, x^2, \ldots y, y^2, \ldots \}$ with $xy=0$, Hesselholt
\cite{Axes} has similarly computed the algebraic $K$--theory of the
coordinate axes $k[x,y]/(xy)$ in terms of the $RO(S^1)$--graded TR--theory of $k$. Computations of representation--graded TR--theory are fundamental to the computation of algebraic $K$--theory, particularly for non-regular schemes.

In the non-equivariant case, computations of TR--groups are possible because the TR--groups form a Witt complex, a rigid algebraic structure. We would like to define a new algebraic
structure, an $R(S^1)$--graded Witt complex, embodying the structure of $R(S^1)$--graded TR. The first step in developing such an algebraic structure is
completely understanding an example of TR$^n_{\alpha}(A;p)$ for
some ring $A$ and all $\alpha \in
R(S^1)$. In this paper we
compute TR$^n_{\alpha}(\mathbb{F}_p;p)$ to provide such an
example.

To study these $R(S^1)$--graded TR--groups, we first consider the
operators and relations that we have in this $R(S^1)$--graded
version of TR (see \cite{KWitt}). Again, inclusion of fixed points induces the
Frobenius map,
\begin{displaymath}
F: \TR^n_{\alpha}(A;p) \rightarrow \TR^{n-1}_{\alpha}(A;p).
\end{displaymath}
As in the integer graded case, this map has an associated
transfer, the Verschiebung
\begin{displaymath}
V: \TR^{n-1}_{\alpha}(A;p) \rightarrow \TR^n_{\alpha}(A;p).
\end{displaymath}

These maps again satisfy the relations $FV=p,$ and  $VF = V(1)$. While these maps F and V look similar to the $\mathbb{Z}$--graded case, the restriction map is quite different. In the $R(S^1)$--graded setting we have
$$
R: \TR^n_{\alpha}(A;p) \rightarrow \TR^{n-1}_{\alpha'}(A;p),
$$
where $\alpha'$ is as defined in the Introduction.  We explain in Section \ref{Spectral} how this restriction map is defined.

For a fixed representation $\alpha$, $\TR^n_{* + \alpha}(A;p)$ is a module over $\TR^n_{*}(A;p)$, where $*$ is an integer grading. In Section \ref{HatGamma} we see that understanding this module structure is essential to our main result.

\section{Norm--Restriction Sequence} \label{Spectral}

\subsection{Construction of norm--restriction diagram}
Throughout this section, we use the notation $T^{-\alpha}$ for $T(A) \wedge S^{-\alpha}$. We are considering $T(A)$ as a genuine $S^1$--spectrum, and all constructions take place in the $S^1$--stable category. Our main tool in doing computations is the fundamental long exact sequence of
TR--groups and its relation to the norm--restriction cofiber sequence.
Consider the cofibration sequence
 \begin{displaymath}
 E_{+} \rightarrow S^0 \rightarrow \tilde{E}
\end{displaymath}
where $E$ denotes a free contractible $S^1$--space and the first
map is given by projection onto the non-basepoint of $S^0$. We can
smash this cofibration sequence with $T^{-\alpha}$ to get
\begin{displaymath}
 E_+\wedge T^{-\alpha}\rightarrow T^{-\alpha}\rightarrow \tilde{E} \wedge T^{-\alpha}.
 \end{displaymath}
 If we then consider the map from $T^{-\alpha}$ to the function spectrum $F(E_+, T^{-\alpha})$ induced by $E_+ \rightarrow S^0$, we get a diagram of cofiber sequences
 \begin{displaymath}
  \xymatrix@C=9pt{
  E_+ \wedge T^{-\alpha}\ar[r] \ar[d] & T^{-\alpha}\ar[r] \ar[d] & \tilde{E} \wedge T^{-\alpha}\ar[d] &  \\
  E_+ \wedge F(E_+, T^{-\alpha}) \ar[r] & F(E_+, T^{-\alpha}) \ar[r]& \tilde{E} \wedge F(E_+, T^{-\alpha}) &\\}
 \end{displaymath}
We can then take the $C_{p^{n-1}}$ fixed point spectra of these spectra
to get a diagram {\small{
\begin{displaymath}
  \xymatrix@C=9pt{
  (E_+ \wedge T^{-\alpha})^{C_{p^{n-1}}}\ar[r] \ar[d] &(T^{-\alpha})^{C_{p^{n-1}}} \ar[r] \ar[d] & (\tilde{E} \wedge T^{-\alpha})^{C_{p^{n-1}}}\ar[d] &\\
   (E_+ \wedge F(E_+, T^{-\alpha}))^{C_{p^{n-1}}}\ar[r] & (F(E_+, T^{-\alpha}))^{C_{p^{n-1}}} \ar[r]& (\tilde{E} \wedge F(E_+, T^{-\alpha}))^{C_{p^{n-1}}}  &\\}
 \end{displaymath} }}
The spectrum in the center of the top line of our diagram
computes $\TR^n_{\alpha}(A ;p)$. We would
also like to identify the top right spectrum as computing an $R(S^1)$--graded
TR--group. The spectrum $T(A)$ is cyclotomic \cite{KWitt}, which means, in particular, that we have an equivalence
of $S^1$--spectra
\begin{displaymath}
r: \rho_p^*(\tilde{E} \wedge T(A))^{C_p}
\rightarrow T(A).
\end{displaymath}
Since $\alpha' = \rho_p^*(\alpha^{C_p})$ this gives an equivalence
$$
\rho_p^*(\tilde{E} \wedge T^{-\alpha})^{C_p}
 \rightarrow T^{-\alpha'}.
$$
Using this, the spectrum
$(\tilde{E} \wedge T^{-\alpha})^{C_{p^{n-1}}}$ computes
$\TR_{\alpha '}^{n-1}(A;p)$. This gives a map
\begin{displaymath}
R: \TR^n_{\alpha}(A;p) \rightarrow \TR^{n-1}_{\alpha
'}(A;p)
\end{displaymath}
in the top row of the diagram. This is the restriction map.
  We now
turn our attention to the bottom row of the diagram. Let
$$
\mathbb{H}^{-q}(C_{p^{n-1}}, T^{-\alpha}) := \pi_q(F(E_+,
T^{-\alpha})^{C_{p^{n-1}}}).
$$
These are the homotopy groups of the group cohomology spectrum of $C_{p^{n-1}}$. This spectrum is also known as the homotopy fixed point spectrum.
Similarly, let
$$
\hat{\mathbb{H}}^{-q}(C_{p^{n-1}},
T^{-\alpha}):= \pi_q((\tilde{E} \wedge F(E_+,
T^{-\alpha}))^{C_{p^{n-1}}}).
$$
 These are the homotopy groups of the
Tate spectrum of $C_{p^{n-1}}$. From \cite[Proposition 2.1]{KWitt},
$$
\pi_q((E_+ \wedge
T^{-\alpha})^{C_{p^{n-1}}}) \cong \pi_q((E_+ \wedge
F(E_+, T^{-\alpha}))^{C_{p^{n-1}}})
$$
and we write
$\mathbb{H}_q(C_{p^{n-1}}, T^{-\alpha}) $ for both of these. Note that if we write $\alpha$ as a formal difference $\alpha = \beta - \gamma$,
$$
\pi_q((E_+ \wedge
T^{-\alpha})^{C_{p^{n-1}}}) \cong [S^q  \wedge S^1/C_{p^{n-1}+} \wedge S^{\beta},  E_+ \wedge T(A)\wedge S^{\gamma}]_{S^1},
$$
hence in the introduction we used the notation $\pi_{\alpha}T(A)_{hC_{p^{n-1}}}$ for the group $\mathbb{H}_0(C_{p^{n-1}}, T^{-\alpha})$.
These
are the homotopy groups of the group homology, or homotopy orbit, spectrum. Returning to our diagram above, we have a diagram of long
exact sequences:
\begin{displaymath}
\xymatrix@C=9pt{
\ar[r] & \mathbb{H}_{0}(C_{p^{n-1}}, T^{-\alpha}) \ar[r] \ar@{=}[d] & \TR^n_{\alpha}(A;p)  \ar[r] \ar[d]^{\Gamma_{n-1}} & \TR^{n-1} _{\alpha '} (A;p)  \ar[r] \ar[d]^{\hat{\Gamma}_n-1} &\\
\ar[r] & \mathbb{H}_{0}(C_{p^{n-1}}, T^{-\alpha}) \ar[r] &
\mathbb{H}^{0}(C_{p^{n-1}}, T^{-\alpha}) \ar[r] &
\hat{\mathbb{H}}^{0}(C_{p^{n-1}}, T^{-\alpha}) \ar[r]&
\\}
\end{displaymath}
The top row of this diagram is the fundamental long exact sequence
of TR--groups. The bottom row is the norm--restriction sequence. We show in Lemma \ref{tate1} below that the bottom row depends only on the
dimension of $\alpha$ and not on the representation itself, and is
thus easier to compute. In particular, the following spectral sequences allow
us to compute the groups on the bottom row:
\begin{displaymath}
\hat{E}^2_{s, t}= \hat{H}^{-s}(C_{p^{n-1}},
\pi_t(T^{-\alpha})) \Rightarrow
\hat{\mathbb{H}}^{-s-t}(C_{p^{n-1}}, T^{-\alpha}),
\end{displaymath}
\begin{displaymath}
{E}^2_{s, t}= {H}_s(C_{p^{n-1}}, \pi_t(T^{-\alpha})) \Rightarrow \mathbb{H}_{s+t}(C_{p^{n-1}},
T^{-\alpha}),
\end{displaymath}
and
\begin{displaymath}
{E}^2_{s, t}= {H}^{-s}(C_{p^{n-1}}, \pi_t(T^{-\alpha})) \Rightarrow \mathbb{H}^{-s-t}(C_{p^{n-1}},
T^{-\alpha}).
\end{displaymath}

\subsection{Computations for $\mathbb{F}_p$}
We now turn our attention to the case $A = \mathbb{F}_p$.
\begin{proposition} \label{diagram}
The norm--restriction diagram of long exact sequences is of one of two forms. If $|\alpha| \geq 0$, the diagram is of the form
$$
\xymatrix{
\ar[r]^-{0}  & \mathbb{Z}/p^{r_n} \ar[r]^-{N} \ar@{=}[d] & \TR^n_{\alpha}(\mathbb{F}_p;p)  \ar[r]^-{R} \ar[d]^{\Gamma_{n-1}} & \TR^{n-1}_{\alpha'}(\mathbb{F}_p;p)  \ar[r] \ar[d]^{\hat{\Gamma}_{n-1}} & \mathbb{Z}/p^{r_n-1}\ar[r] \ar@{=}[d] &\\
\ar[r]^-{0} & \mathbb{Z}/p^{r_n} \ar[r]  & \mathbb{Z}/p^n \ar[r]
& \mathbb{Z}/p^{n-1}  \ar[r]& \mathbb{Z}/p^{r_n-1} \ar[r]&\\}
$$
where $r_n=\min(n, |\alpha| + 1)$. If $|\alpha|<0$, the diagram is of the form
$$
\xymatrix{
\ar[r]  & 0 \ar[r]^-{N} \ar@{=}[d] & \TR^n_{\alpha}(\mathbb{F}_p;p)  \ar[r]^-{R} \ar[d]^{\Gamma_{n-1}} & \TR^{n-1}_{\alpha'}(\mathbb{F}_p;p)  \ar[r] \ar[d]^{\hat{\Gamma}_{n-1}} & 0 \ar[r] \ar@{=}[d] &\\
\ar[r] & 0 \ar[r]  & \mathbb{Z}/p^{n-1} \ar[r]
& \mathbb{Z}/p^{n-1}  \ar[r]& 0 \ar[r]&\\}
$$
\end{proposition}

\begin{corollary}\label{NegR}
If $|\alpha|<0$,
$$
\TR^n_{\alpha}(\mathbb{F}_p;p)  \cong \TR^{n-1}_{ \alpha'}(\mathbb{F}_p;p).
$$
\end{corollary}

We prove this proposition by using the above spectral sequences to evaluate the homotopy groups of the Tate spectrum, the group homology spectrum, and the group cohomology spectrum. We start with the Tate spectrum.
\subsubsection{The Tate spectrum}

\begin{lem} \label{tate1}
Letting $[S^{-\alpha}]$ denote a shift in degree by $-2|\alpha|$,
\begin{displaymath}
\hat{\mathbb{H}}^*(C_{p^{n-1}}, T(\mathbb{F}_p)\wedge S^{-\alpha})
\cong S_{\mathbb{Z}/p^{n-1}}\{\hat{\sigma},
\hat{\sigma}^{-1}\}[S^{-\alpha}], |\hat{\sigma}| = 2.
\end{displaymath}
\end{lem}
To prove this lemma, we recall some results of Hesselholt and Madsen. Hesselholt and Madsen
\cite[Section 5]{KWitt}  have computed that the spectral sequence
\begin{displaymath}
\hat{E}^2_{s, t}= \hat{H}^{-s}(C_{p^{n-1}}, \pi_t(T(\mathbb{F}_p)))
\Rightarrow \hat{\mathbb{H}}^{-s-t}(C_{p^{n-1}}, T(\mathbb{F}_p))
\end{displaymath}
 has $\hat{E}^2$ term
 \begin{displaymath}
\hat{E}^2 = \Lambda_{\mathbb{F}_p}\{u_{n-1}\} \otimes
S_{\mathbb{F}_p}\{t, t^{-1}\} \otimes S_{\mathbb{F}_p}\{\sigma \}
 \end{displaymath}
 where the classes $u_{n-1}, t, \sigma$ are in bidegrees $ (-1, 0), (-2, 0), (0,2)$ respectively. The notation $S_{\mathbb{F}_p}\{t, t^{-1} \} $ denotes a polynomial algebra over $\mathbb{F}_p$ with generators $t$ and $t^{-1}$. Hesselholt and Madsen also computed that the non-zero differentials are given by
 $$
 d^{2n-1}u_{n-1} = t^{n}\sigma^{n-1}
 $$
and the extensions in the passage from $\hat{E}^{\infty}$ to the homotopy groups are maximally nontrivial. Hence the spectral sequence converges to
\begin{displaymath}
\hat{\mathbb{H}}^*(C_{p^{n-1}}, T(\mathbb{F}_p)) \cong
S_{\mathbb{Z}/p^{n-1}}\{\hat{\sigma}, \hat{\sigma}^{-1}\},
|\hat{\sigma}|=2.
\end{displaymath}
In our case we would like to look at the spectral sequence
\begin{displaymath}
\hat{E}^2_{s, t}= \hat{H}^{-s}(C_{p^{n-1}},
\pi_t(T^{-\alpha})) \Rightarrow
\hat{\mathbb{H}}^{-s-t}(C_{p^{n-1}}, T^{-\alpha}).
\end{displaymath}
We now verify that the $E^2$--term and differentials in this spectral sequence depend only on the dimension of $\alpha$, and not on which virtual representation it is of that dimension. We first observe that we can write the $E^2$--term as
$$
\hat{E}^2_{s, t}= \hat{H}^{-s}(C_{p^{n-1}},
\pi_t(T(\mathbb{F}_p)\wedge S^{-\alpha})) \cong \hat{H}^{-s}(C_{p^{n-1}},
\pi_{t+2|\alpha|}(T(\mathbb{F}_p))).
$$
Since the $C_{p^{n-1}}$--action on $T(\mathbb{F}_p)$ comes from an $S^1$ action, it is trivial on homotopy, and hence this $E^2$--term doesn't depend on any information about the representation $\alpha$ other than its dimension. Thus we write
$$
\hat{E}^2(T(\mathbb{F}_p) \wedge S^{-\alpha}) = \hat{E}^2(T(\mathbb{F}_p))[S^{-\alpha}] = \big(\Lambda_{\mathbb{F}_p}\{u_{n-1}\} \otimes
S_{\mathbb{F}_p}\{t, t^{-1}\} \otimes S_{\mathbb{F}_p}\{\sigma \}\big)[S^{-\alpha}]
$$
where $[S^{-\alpha}]$ denotes a shift in bidegree by $(0, -2|\alpha|)$.
We now look at the differentials in this spectral sequence.
\begin{lem}
 For $i, j\in \mathbb{Z}, j\geq0$, the classes $t^i\sigma^j[S^{-\alpha}] \in \hat{E}^2(C_{p^{n-1}}, T^{-\alpha})$ are permanent cycles.
\begin{proof}
By \cite[Section 4]{local},
$$
\hat{E}^2(S^1, T(\mathbb{F}_p)) \cong S_{\mathbb{F}_p}\{t, t^{-1}\} \otimes S_{\mathbb{F}_p}({\sigma}).
$$
Thus,
$$
\hat{E}^2(S^1, T(\mathbb{F}_p)\wedge S^{-\alpha}) = \hat{E}^2(S^1, T(\mathbb{F}_p))[S^{-\alpha}]
$$
is also concentrated in even total degree. Hence there are no differentials in this spectral sequence. The canonical inclusion induces a map of spectral sequences
$$
\hat{E}^2(S^1, T^{-\alpha}) \rightarrow \hat{E}^2(C_{p^{n-1}}, T^{-\alpha}).
$$
By \cite[Section 4]{local} this map has image
$$
\big(S_{\mathbb{F}_p}\{t, t^{-1}\} \otimes S_{\mathbb{F}_p}({\sigma})\big)[S^{-\alpha}] \in \hat{E}^2(C_{p^{n-1}}, T^{-\alpha})
$$
and hence these classes are permanent cycles.
\end{proof}
\end{lem}
So we conclude that either
$$
d^r(u_{n-1} t^i\sigma^j[S^{-\alpha}]) \neq 0
$$
for some $r \geq 2$, or all differentials are zero. Suppose for some $r, i, j$ this differential is nonzero. Note that
$$
\hat{E}^2_{s,t}(C_{p^{n-1}}, T(\mathbb{F}_p) \wedge S^{-\alpha}) = \hat{H}^{-s}(C_{p^{n-1}},
\pi_{t+2|\alpha|}(T(\mathbb{F}_p))) = 0
$$
if $t$ is odd. Hence $r$ must be odd. So, we rewrite our non-zero differential as
$$
d^{2r+1}(u_{n-1} t^i\sigma^j[S^{-\alpha}]) = d^{2r+1}(u_{n-1})t^i\sigma^j[S^{-\alpha}].
$$
Hesselholt and Madsen prove in \cite[Lemma 5.4]{KWitt}  that
$$
d^{2r+1}(u_{n-1}) \neq 0
$$
only when $r=n-1$ and in that case
$$
d^{2n-1}(u_{n-1}) = t^{n}\sigma^{n-1}.
$$
Thus the only nonzero differentials in the spectral sequence we are studying are generated by
$$
d^{2n-1}(u_{n-1} [S^{-\alpha}]) = t^{n}\sigma^{n-1}[S^{-\alpha}].
$$
The extensions in the passage from $\hat{E}^{\infty}$ are maximally nontrivial, so this spectral sequence converges to
\begin{displaymath}
\hat{\mathbb{H}}^*(C_{p^{n-1}}, T(\mathbb{F}_p)\wedge S^{-\alpha})
\cong S_{\mathbb{Z}/p^{n-1}}\{\hat{\sigma},
\hat{\sigma}^{-1}\}[S^{-\alpha}], |\hat{\sigma}| = 2.
\end{displaymath}
This proves Lemma \ref{tate1}.
\subsubsection{The group homology spectrum}

We now turn our attention to the group homology spectrum.
\begin{lem} \label{homology}
The homotopy groups of the group homology spectrum are given by
$$
\mathbb{H}_{0}(C_{p^{n-1}}, T^{-\alpha}) \cong \left\{ \begin{array}{ll} \mathbb{Z}/p^{r_n} & \textup{  if  }  |\alpha| \geq 0 \\
 0 & \textup{  if } |\alpha| <0 \\ \end{array} \right.
$$
and
$$
\mathbb{H}_{-1}(C_{p^{n-1}}, T^{-\alpha}) \cong \left\{ \begin{array}{ll} \mathbb{Z}/p^{r_n-1} & \textup{  if  } |\alpha| \geq 0 \\
 0 & \textup{  if } |\alpha| <0 \\ \end{array} \right.
$$
where $r_n=\min(n,  |\alpha| + 1)$.
\begin{proof}
Recall from above
that we have a spectral sequence
\begin{displaymath}
{E}^2_{s, t}= {H}_s(C_{p^{n-1}}, \pi_t(T(\mathbb{F}_p)\wedge
S^{-\alpha})) \Rightarrow \mathbb{H}_{s+t}(C_{p^{n-1}},
T(\mathbb{F}_p)\wedge S^{-\alpha}).
\end{displaymath}
We can see easily from the spectral sequence for $\mathbb{H}_{0}$ (and similarly for $\mathbb{H}_{-1}$)
that
$$
\mathbb{H}_{0}(C_{p^{n-1}}, T(\mathbb{F}_p) \wedge S^{-\alpha}) = 0
\textup{   for   }  |\alpha| <0.
$$
We now compute these groups for $ |\alpha|\geq 0$. Note that
$$
E^2_{s,t} =  {H}_s(C_{p^{n-1}}, \pi_t(T(\mathbb{F}_p)\wedge
S^{-\alpha})) \cong
{H}_s(C_{p^{n-1}},
\pi_{t+2|\alpha|}(T(\mathbb{F}_p)).
$$
B\"{o}kstedt and Breen (\cite{Bok2}, \cite{Breen}) calculated that $\pi _*(T(\mathbb{F}_p)) \cong
S_{\mathbb{F}_p}\{\sigma_1\}, |\sigma_1|=2$. So our $E^2$ term is
$$
E^2_{s,t} = \left\{ \begin{array}{ccccc}
H_s(C_{p^{n-1}}, \mathbb{Z}/p) & = & \mathbb{Z}/p & \textup{if} & t+2|\alpha| \geq 0, \textup{even}, s\geq 0 \\
0 & & \textup{else} &
\end{array} \right.
$$
Viewing this $E^2$ term as the shifted first quadrant part of the
$\hat{E}^2$ term of the spectral sequence computing the Tate
cohomology, we have names for these classes. We have
$$
E^2_{s,t}(C_{p^{n}}, T(\mathbb{F}_p) \wedge S^{-\alpha}) =
\big(S_{\mathbb{F}_p}\{{t^{-1}}\} \otimes S_{\mathbb{F}_p}\{{\sigma}\}
\otimes \Lambda_{\mathbb{F}_p}({t^{-1}u_n})\big)[S^{-\alpha}].
$$
In our discussion of the Tate spectrum above, we saw that the differentials in
the spectral sequence are generated by
$$
d^{2n+1}u_n =
t^{n+1}\sigma ^n.
$$ Using this and the fact that the extensions in
the passage from $E^\infty$ to homotopy groups are maximally
nontrivial, the stated result follows.
\end{proof}
\end{lem}
\subsubsection{The group cohomology spectrum}

Finally we look at the homotopy groups of the group cohomology
spectrum.
\begin{lem} \label{cohomology}
The homotopy groups of the group cohomology spectrum are given by
$$
\mathbb{H}^{0}(C_{p^{n-1}}, T^{-\alpha}) =
\left\{\begin{array}{cc}
\mathbb{Z}/p^{n} &  |\alpha|\geq 0   \\
\mathbb{Z}/p^{n-1} & |\alpha| <0  \\
 0 & \textup{else}
\end{array} \right.
$$
\begin{proof}
Recall that we have a spectral sequence
\begin{displaymath}
{E}^2_{s, t}= {H}^{-s}(C_{p^{n-1}}, \pi_t(T^{-\alpha})) \Rightarrow \mathbb{H}^{-s-t}(C_{p^{n-1}},
T^{-\alpha}).
\end{displaymath}
If we look at the spectral sequence
\begin{displaymath}
{E}^2_{s, t}= {H}^{-s}(C_{p^{n-1}}, \pi_t(T(\mathbb{F}_p)))
\Rightarrow \mathbb{H}^{-s-t}(C_{p^{n-1}}, T(\mathbb{F}_p))
\end{displaymath}
we see that the $E^2$--term of this spectral sequence is the second
quadrant part of the $\hat{E}^2$--term for our Tate spectral
sequence. In other words, the $E^2$--term for this group cohomology
spectral sequence is
 $$ E^2 = \Lambda_{\mathbb{F}_p}(u_{n-1}) \otimes S_{\mathbb{F}_p}\{t\} \otimes
S_{\mathbb{F}_p}\{\sigma\}.$$ Again, the differentials in this
spectral sequence are generated by $d^{2n-1}u_{n-1} =
t^{n}\sigma^{n-1}.$ Hence we have
$$
E^{\infty} =S_{\mathbb{F}_p}\{t\}
\otimes S_{\mathbb{F}_p}\{\sigma\}/t^{n}\sigma^{n-1}.
$$
As above, in our case this is shifted by twice the complex dimension of our representation
$\alpha$. The desired result follows by \cite[Section 5]{KWitt}.
\end{proof}
\end{lem}

The results of Lemmas \ref{tate1}, \ref{homology}, and \ref{cohomology}, prove Proposition \ref{diagram}. Throughout our computations we will carefully study the map
$$\hat{\Gamma}_n: \TR^{n}_{\alpha}(\mathbb{F}_p;p) \rightarrow
\hat{\mathbb{H}}^{0}(C_{p^{n}}, T(\mathbb{F}_p) \wedge
S^{-\beta})
$$
for $\beta \in R(S^1)$ with $\beta' = \alpha$. Hesselholt and Madsen have proven the following proposition about this map \cite[Addendum 9.1]{KWitt}.

\begin{proposition}\label{GammaHat}
The map $\hat{\Gamma}_n$ is an isomorphism if the integers $ |\alpha|, |\alpha'|,
|\alpha''|, \ldots |\alpha^{(n-1)}|$ are all nonnegative.
\end{proposition}

Note that in the cited reference the proposition is stated in the particular case of representations of the form $\alpha = q - \lambda$, $q \in \mathbb{Z}$, $\lambda$ an actual complex representation of $S^1$. The proof extends to the general case of a virtual representation $\alpha \in R(S^1)$, in which case the statement of the result is as above. Since for an actual representation $\lambda$, $ |\lambda| \geq|\lambda'| \geq \ldots |\lambda^{(n-1)}|$, the proposition in the cited reference is stated as: $\hat{\Gamma}_n$ is an isomorphism if $q \geq |\lambda|$.

\section{Inductive method}\label{compute}

The main result of this paper is the computation of $\TR^n_{\alpha}(\mathbb{F}_p;p)$, for all $n \in \mathbb{N}$ and $\alpha \in R(S^1)$. We prove the main theorem inductively in two parts. Given the group $\TR^{n-1}_{\alpha'}(\mathbb{F}_p;p)$ and the map $\hat{\Gamma}_{n-1}$, we compute the group $\TR^n_{\alpha}(\mathbb{F}_p;p)$ and the map $\hat{\Gamma}_n$. We begin with the base case for our induction.

\subsection{Base case}

The base case of our induction is the computation of $\TR^1_{\alpha^{(n-1)}}(\mathbb{F}_p;p)$. Note that
$$
\TR^1_{\alpha^{(n-1)}}(\mathbb{F}_p;p) \cong
\pi_{2|\alpha^{(n-1)}|}(T(\mathbb{F}_p)) \cong \left\{ \begin{array}{cc}
 \mathbb{Z}/p& |\alpha^{(n-1)}| \geq 0  \\
 0 & \textup{  else  }
 \end{array}\right.
$$
This is a shift of the homotopy groups
$$
\TR^1_*(\mathbb{F}_p; p) \cong \pi_*(T(\mathbb{F}_p)) \cong
S_{\mathbb{Z}/p}\{\sigma_1\}
$$
where $|\sigma_1| = 2$. By Proposition \ref{GammaHat} the map $\hat{\Gamma}_1$ is injective.

\subsection{Inductive step}

We now compute the group $\TR^n_{\alpha}(\mathbb{F}_p;p)$ from the group $\TR^{n-1}_{\alpha'}(\mathbb{F}_p;p)$ and the map $\hat{\Gamma}_{n-1}$.

\begin{proposition}\label{Pullback1}

Let $r_n = \min(n, |\alpha| + 1)$. If $r_n \geq 1$, the group $\TR^n_{\alpha}(\mathbb{F}_p;p)$ is given by the following pullback:
$$
\xymatrix{
\TR^n_{\alpha}(\mathbb{F}_p;p) \ar[r]^{R} \ar[d]^{{\Gamma}_{n-1}} &\TR^{n-1}_{\alpha'}(\mathbb{F}_p;p) \ar[d]^{\hat{\Gamma}_{n-1}} \\
\mathbb{Z}/p^n \ar[r]^{R^h_n} & \mathbb{Z}/p^{n-1} \\}
$$
where $R^h_n(1) \doteq p^{r_n - 1}$. If $r_n <1$, $\TR^n_{ \alpha}(\mathbb{F}_p;p) \cong \TR^{n-1}_{ \alpha'}(\mathbb{F}_p;p) $.

\begin{proof}
By Corollary \ref{NegR}, if $r_n <1$, $\TR^n_{ \alpha}(\mathbb{F}_p;p) \cong \TR^{n-1}_{ \alpha'}(\mathbb{F}_p;p) $. We consider the case where $r_n \geq 1$. We first recall some general facts about diagrams of exact sequences. Suppose we have a diagram of the form
$$
\xymatrix{
0 \ar[r] & A  \ar@{=}[d] & & C  \ar[r] \ar[d] & D  \ar[r] \ar@{=}[d] &\\
0 \ar[r] & A' \ar[r] &  B' \ar[r]  & C' \ar[r]
& D'  \ar[r]& \\}
$$
The pullback
$$
\xymatrix{
R \ar[r] \ar[d] & C \ar[d] &\\
B' \ar[r] & C' &\\
}
$$
makes the above diagram commute, and the sequence
$$
\xymatrix{ 0 \ar[r] & A \ar[r] & R \ar[r] & C \ar[r] & D \ar[r] &\\}
$$
is exact. Any group $B$ making the diagram of exact sequences commute factors through the pullback, giving a diagram of exact sequences
$$
\xymatrix{
0 \ar[r] & A \ar[r] \ar@{=}[d] & B \ar[r] \ar[d] & C  \ar[r] \ar@{=}[d] & D  \ar[r] \ar@{=}[d] &\\
0 \ar[r] & A \ar[r] &  R \ar[r]  & C \ar[r]
& D  \ar[r]& \\}
$$
Then, by the five lemma, $R \cong B$. So, up to isomorphism there is only one group fitting into such a diagram of long exact sequences.

To calculate $\TR^n_{\alpha}(\mathbb{F}_p; p)$ we use the norm--restriction diagram from Section \ref{Spectral}. By Proposition \ref{diagram}, when $r_n \geq 1$, this diagram has the form
$$
\xymatrix{
\ar[r]^-{0}  & \mathbb{Z}/p^{r_n} \ar[r]^-{N} \ar@{=}[d] & \TR^n_{\alpha}(\mathbb{F}_p;p)  \ar[r]^-{R} \ar[d]^{\Gamma_{n-1}} & \TR^{n-1}_{\alpha'}(\mathbb{F}_p;p)  \ar[r]^-{\partial} \ar[d]^{\hat{\Gamma}_{n-1}} & \mathbb{Z}/p^{r_n-1}\ar[r] \ar@{=}[d] &\\
\ar[r]^-{0} & \mathbb{Z}/p^{r_n} \ar[r]^-{N^h}  & \mathbb{Z}/p^n \ar[r]^-{R^h}
& \mathbb{Z}/p^{n-1}  \ar[r]^-{\partial^h}& \mathbb{Z}/p^{r_n-1} \ar[r]&\\}
$$

The maps in the bottom row are given by $N^h(1) \doteq p^{n-r_n}$, $R^h(1) \doteq p^{r_n-1}$ and $\partial^h(1) \doteq 1$. By the arguments above, the group $\TR^n_{\alpha}(\mathbb{F}_p;p)$ is given by the pullback of the diagram
$$
\xymatrix{ & \TR^{n-1}_{\alpha'}(\mathbb{F}_p;p)  \ar[d]^{\hat{\Gamma}_{n-1}} & \\
\mathbb{Z}/p^n \ar[r]^-{R^h}
& \mathbb{Z}/p^{n-1}&
\\}
$$
This completes the proof.
\end{proof}
\end{proposition}

Before we proceed to compute $\hat{\Gamma}_n$, it will help to explicitly compute the pullback diagram above. The following proposition evaluates this pullback.
\begin{proposition}\label{pullback}

Let  $\hat{\Gamma}_{i,n-1}$ be a map of abelian groups
$$
\hat{\Gamma}_{i,n-1}: \mathbb{Z}/p^{l_{i,n-1}} \rightarrow \mathbb{Z}/p^{n-1}.
$$
Suppose $k_{i,n-1} = \textup{length}_{\mathbb{Z}_p}(\textup{ker }\hat{\Gamma}_{i,n-1})$, and $k_{1,n-1} \leq k_{2, n-1} \leq \ldots \leq k_{n-1,n-1}$. Then the following square is a pullback diagram
$$
\xymatrix{
\mathbb{Z}/p^{l_{1, n}} \oplus \mathbb{Z}/p^{l_{2,n}} \oplus \cdots \oplus \mathbb{Z}/p^{l_{n,n}}  \ar[r]^-{R} \ar[d]^{{\Gamma}_{1,n-1} + \ldots + {\Gamma}_{n-1, n-1}} &  \mathbb{Z}/p^{l_{1, n-1}} \oplus \mathbb{Z}/p^{l_{2,n-1}} \oplus \cdots \oplus \mathbb{Z}/p^{l_{n-1,n-1}}  \ar[d]^{\hat{\Gamma}_{1,n-1} + \ldots + \hat{\Gamma}_{n-1, n-1}} &\\
\mathbb{Z}/p^n \ar[r]^{R^h} & \mathbb{Z}/p^{n-1}&\\
}
$$
where the map $R^h$ is given by $R^h(1) \doteq p^{r_n-1}$ and the integers $l_{i,n},$ $ 1 \leq i \leq n$ are defined as follows:
Let  $g_{i,n-1} = l_{i,n-1} - k_{i,n-1}$.  We choose a permutation $\tau \in \Sigma_{n-1}$ satisfying the conditions
$$
r_n \leq k_{\tau(1), n-1} \leq k_{\tau(2), n-1} \leq \ldots k_{\tau(w), n-1} ,
$$
$$
r _n> k_{\tau(w+1), n-1}, k_{\tau(w+2), n-1} , \ldots k_{\tau(n-1), n-1} ,
$$
$$
g_{\tau(w+1),n-1} \geq g_{\tau(w+2), n-1} \geq \ldots \geq g_{\tau(n-1), n-1},
$$
Then
$$
l_{1,n} = \left\{ \begin{array}{ll}r_n & \textup{  if  } w=n-1\\
 \min(n, r_n+ g_{\tau(w+1), n-1}) & \textup{ else}  \\ \end{array} \right.
$$
For $1\leq m \leq w$, $l_{m+1, n} = l_{\tau(m), n-1}$. For $w+2 \leq m \leq n-1$:
$$
l_{m,n}= \left\{ \begin{array}{ll} l_{\tau(m), n-1} & \textup{  if } \tau(v) < \tau(m) \textup{ for some } w+1 \leq v<m  \\
g_{\tau(m),n-1} + k_{\tau(z),n-1} & \textup{ else} \\ \end{array}\right.
$$
Here $\tau(z)$ is the smallest number such that $w+1 \leq z <m$. Finally, $l_{n,n} = k_{1,n-1}$ if $w<n-1$.
\begin{proof}
Let $G$ denote the pullback of the above diagram. We would like to prove that
$$
G \cong \mathbb{Z}/p^{l_{1, n}} \oplus \mathbb{Z}/p^{l_{2,n}} \oplus \cdots \oplus \mathbb{Z}/p^{l_{n,n}}.
$$
By definition
$$
G \cong (x, y_1, y_2, \ldots y_{n-1}) \in \mathbb{Z}/p^n \oplus \mathbb{Z}/p^{l_{1, n-1}} \oplus \cdots \oplus \mathbb{Z}/p^{l_{n-1,n-1}}
$$
such that
$$
p^{r_n-1}x = p^{n-1-g_{1,n-1}}y_1 + p^{n-1-g_{2,n-1}}y_2 + \ldots + p^{n-1-g_{n-1,n-1}}y_{n-1}.
$$
We choose a permutation $\tau \in \Sigma_{n-1}$ as in the statement of the proposition. For each $1 \leq m  \leq w$ the element
$$
\zeta_{\tau(m)} = (p^{n-g_{\tau(m),n-1} - r_n}, 0, \ldots 1, 0,\ldots 0) \in  \mathbb{Z}/p^n \oplus \mathbb{Z}/p^{l_{1,n-1}} \oplus \cdots \oplus \mathbb{Z}/p^{l_{n-1,n-1}}
$$
where the 1 is in the $y_{\tau(m)}$--coordinate, generates a subgroup of $G$ of order $p^{l_{\tau(m),n-1}}.$ Note that the different $m$ give linearly independent elements of $G$. The element
$$
(p^{n - r_n}, 0, \ldots 0) \in  \mathbb{Z}/p^n \oplus \mathbb{Z}/p^{l_{1, n-1}} \oplus \cdots \oplus \mathbb{Z}/p^{l_{n-1,n-1}}
$$
generates a subgroup of $G$ of order $p^r$. If $w=n-1$,  $k_{i, n-1} \geq r_n$ for all $1 \leq i \leq n-1$, so this element is linearly independent from the elements $\zeta_i$. So in the case $w=n-1$ we have computed
$$
G \cong \mathbb{Z}/p^r \oplus \mathbb{Z}/p^{l_{1,n-1}} \oplus \cdots \oplus \mathbb{Z}/p^{l_{n-1,n-1}}.
$$
Suppose then that $w \neq n-1$. For $w+1 \leq m \leq n-1$, if $k_{\tau(m), n-1} = l_{\tau(m),n-1}$, then the element
$$
\eta_{\tau(m)}= (0, 0 \ldots 1, 0 \ldots 0) \in  \mathbb{Z}/p^n \oplus \mathbb{Z}/p^{l_{1,n-1}} \oplus \cdots \oplus \mathbb{Z}/p^{l_{n-1,n-1}}
$$
which is non-zero only in the $y_{\tau(m)}$--coordinate, generates a subgroup of length $l_{\tau(m),n-1}$. For $k_{\tau(m),n-1} < l_{\tau(m),n-1}$, $w+1 \leq m \leq n-1$, if $l_{\tau(m),n-1} -k_{\tau(m),n-1} + r_n \leq n$, we consider  the element
$$
\psi_{\tau(m)} = (p^{n-g_{\tau(m),n-1} - r_n}, 0, \ldots 1, 0,\ldots 0) \in  \mathbb{Z}/p^n \oplus \mathbb{Z}/p^{l_{1,n-1}} \oplus \cdots \oplus \mathbb{Z}/p^{l_{n-1,n-1}}
$$
where the 1 is in the $y_{\tau(m)}$--coordinate. This generates an order $p^{r+g_{\tau(m), n-1}}$ subgroup of $G$. If instead  $g_{\tau(m),n-1} + r_n \geq n$, we consider the element
$$
\psi_{\tau(m)} = (1, 0, \ldots p^{g_{\tau(m),n-1} +r_n - n}, 0 \ldots 0) \in  \mathbb{Z}/p^n \oplus \mathbb{Z}/p^{l_{1,n-1}} \oplus \cdots \oplus \mathbb{Z}/p^{l_{n-1,n-1}}
$$
where the second nonzero entry is in the $y_{\tau(m)}$--coordinate. This generates a subgroup of length $n$. For $w+2 \leq m \leq n-1$, $w+1 \leq m' <m$, consider also the elements
$$
\nu_{\tau(m'), \tau(m)} = (0, 0, \ldots, p^{g_{\tau(m'), n-1}  - g_{\tau(m), n-1 }}, 0 \ldots  -1,0, \ldots 0).
$$
Here, the first nonzero term is in the $y_{\tau(m')}$--coordinate, and the second nonzero term is in the $y_{\tau(m)}$--coordinate. These generate subgroups of length $\max(l_{\tau(m), n-1}, g_{\tau(m), n-1} + k_{\tau(m'),n-1})$. Since we have ordered the summands such that $k_{\tau(m),n-1} \leq k_{\tau(m'),n-1}$ exactly when $\tau(m) \leq \tau(m')$, we conclude that the length of this group is $l_{\tau(m),n-1}$ if $\tau(m') \leq \tau(m)$, and $g_{\tau(m),n-1} + k_{\tau(m'),n-1}$ otherwise.
Finally, consider the element
$$
\nu_{n-1} = (0, p^{g_{1,n-1}}, 0, 0 \ldots 0).
$$
which generates a subgroup of order $p^{k_{1,n-1}}$.

Together, the elements above generate the group $G.$ In order to identify this group, we need to find a linearly independent set of generators that spans.  For a fixed $w+2 \leq m \leq n-1$ let $m'$ be given by $\tau(m') = \min(\tau(i) | w+1 \leq i <m)$. Our set of generators is the following:
$$
\{\psi_{\tau(w+1)} \}, \{\zeta_{\tau(m)}| 1 \leq m  \leq w\},
$$
$$
\{\eta_{\tau(m)} | w+2 \leq m \leq n-1,  k_{\tau(m),n-1} \geq l_{\tau(m),n-1}\},
$$
$$
\{\nu_{\tau(m), \tau(m')} | w+2 \leq m \leq n-1,  k_{\tau(m),n-1} < l_{\tau(m),n-1}\}, \{\nu_{n-1}\}.
$$
The summands in the proposition correspond to these generators in the following way:
$$
\psi_{\tau(w+1)} \mapsto \mathbb{Z}/p^{l_{1,n}}, \textup{        }\zeta_{\tau(m)} \mapsto \mathbb{Z}/p^{l_{m+1, n}}, \textup{         } \eta_{\tau(m)} \mapsto \mathbb{Z}/p^{l_{m, n}},
$$
$$
\nu_{\tau(m), \tau(m')} \mapsto \mathbb{Z}/p^{l_{m,n}}, \textup{         }\nu_{n-1} \mapsto \mathbb{Z}/p^{l_{n, n}}.
$$
Then the values of the lengths $l_{i,n}$ follow directly from our analysis above.
\end{proof}
\end{proposition}
We can choose an isomorphism
$$
\iota_n: \mathbb{Z}/p^{l_{1,n}} \oplus \cdots \oplus \mathbb{Z}/p^{l_{n,n}} \rightarrow \TR^n_{\alpha}(\mathbb{F}_p;p).
$$
The proof of Proposition \ref{pullback} above determines the restriction map:
$$
\xymatrix{\mathbb{Z}/p^{l_{1,n}} \oplus \cdots \oplus \mathbb{Z}/p^{l_{n,n}} \ar[rr]^-{\iota_n}_-{\sim} \ar[d]^{R^{\iota}}&& \TR^n_{\alpha}(\mathbb{F}_p;p) \ar[d]^{R} & \\
\mathbb{Z}/p^{l_{1,n-1}} \oplus \cdots \oplus \mathbb{Z}/p^{l_{n-1,n-1}} \ar[rr]^-{\iota_{n-1}}_-{\sim} && \TR^{n-1}_{\alpha'}(\mathbb{F}_p;p)& \\}
$$
Note that we can write the isomorphism $\iota_n$ as
$$
\iota_n = \iota_{1, n} + \ldots + \iota_{n,n}
$$
where
$$
\iota_{j,n}: \mathbb{Z}/p^{l_{j,n}} \rightarrow \TR^n_{\alpha}(\mathbb{F}_p;p).
$$
Then we can express the composite $R \iota_{j,n}$ as
$$
R\iota_{j,n}(a) = \sum_{1 \leq i \leq n-1} \iota_{i,n-1} (R^{\iota}_{i,j,n}(a))
$$
where
$$
R^{\iota}_{i,j,n}: \mathbb{Z}/p^{l_{j,n}} \rightarrow \mathbb{Z}/p^{l_{i,n-1}}.
$$
We now determine the map $R$ by specifying the maps $R^{\iota}_{i,j,n}$.
\begin{thm}\label{restriction}
The maps $R^{\iota}_{i,j,n}$ are given by the following. If $i=\tau(w+1)$:
$$
R^{\iota}_{i, 1, n}(1) \doteq \left\{ \begin{array}{ll}
1 & \textup{ if  }  r_n + g_{\tau(w+1), n-1}  \leq n  \\
p^{r_n + g_{\tau(w+1),n-1}  - n} & \textup{ if  }r_n + g_{\tau(w+1), n-1}  > n \\ \end{array} \right.
$$
For $i \neq \tau(w+1)$,  $R^{\iota}_{i, 1, n}(1) = 0$. For $1 < m \leq n$, if $m-1 \leq w$:
$$
R^{\iota}_{i, m, n}(1) \doteq \left\{ \begin{array}{ll} 1  & \textup{ if  } i=\tau(m-1)  \\
 0 & \textup{ else  }  \\ \end{array} \right.
$$
For $w<m -1<n$:
$$
R^{\iota}_{i, m, n}(1) \doteq \left\{ \begin{array}{ll} 1 & \textup{ if  }i=\tau(m),  k_{\tau(m), n-1} = l_{\tau(m), n-1}  \\
-1 & \textup{ if  }i=\tau(m), k_{\tau(m), n-1} < l_{\tau(m),n-1}   \\
p^{ g_{\tau(m'),n-1} - g_{\tau(m), n-1} } & \textup{ if  } i=\tau(m'), k_{\tau(m), n-1} < l_{\tau(m),n-1} \\
0 &  \textup{ else  } \\ \end{array} \right.
$$
If $w <n-1$:
$$
R^{\iota}_{i, n, n}(1) \doteq \left\{ \begin{array}{ll} p^{g_{1, n-1}}  & \textup{ if  } i=1  \\
 0 & \textup{ else  }  \\ \end{array} \right.
$$
\begin{proof}
This follows directly from the proof of Proposition \ref{pullback}.
\end{proof}
\end{thm}
The proof of Proposition \ref{pullback} also determines the map
$$
\xymatrix{\Gamma_{n-1}: & \mathbb{Z}/p^{l_{1, n}} \oplus \mathbb{Z}/p^{l_{2,n}} \oplus \cdots \oplus \mathbb{Z}/p^{l_{n,n}}  \ar[r]^-{\iota_n}_-{\sim} & \TR^n_{\alpha}(\mathbb{F}_p;p) \ar[r]^-{\Gamma_{n-1}} & \mathbb{Z}/p^{n} &\\}
$$
\begin{proposition}\label{Gamma}
Let $\overline{k}=\textup{length}_{\mathbb{Z}_p}(\textup{ ker } {\Gamma}_{i,n-1}).$ The maps $\Gamma_{i, n-1}: \mathbb{Z}/p^{l_{i, n}} \rightarrow \mathbb{Z}/p^{n}$ are determined up to units by the following:
$$
\overline{k}_{i,n} = \left\{ \begin{array}{lll} 0 & \textup{if} & i=1 \\
k_{\tau(i-1), n-1} - r_n & \textup{if} & k_{\tau(i-1), n-1} \geq r_n \\
 l_{i,n} & \textup{else}
\end{array}\right.
$$
\begin{proof}
This follows directly from the proof of Proposition \ref{pullback}
\end{proof}
\end{proposition}

Given the group $\TR^{n-1}_{\alpha'}(\mathbb{F}_p;p)$ and the map $\hat{\Gamma}_{n-1}$ on this group, we have now calculated the group $\TR^n_{\alpha}(\mathbb{F}_p;p)$. In order to complete the induction step we also need to compute the map $\hat{\Gamma}_n$. We study this map in the next section.

\section{The map $\hat{\Gamma}_n$}\label{HatGamma}

Proposition \ref{GammaHat} tells us how the map $\hat{\Gamma}_n$ behaves for virtual representations $\alpha$ where all of  $|\alpha|, |\alpha'|, |\alpha''|, \ldots |\alpha^{(n-1)}|$ are nonnegative. To understand the map $\hat{\Gamma}_n$ when at least one of $|\alpha|, |\alpha'|, |\alpha''|, \ldots |\alpha^{(n-1)}|$ is negative, we study the module structure of
$\TR^{n}_{*+\alpha}(\mathbb{F}_p;p)$ as a
$\TR^{n}_*(\mathbb{F}_p; p)$--module where $*$ denotes an even integer grading. We
first recall from ~\cite[Proposition 5.4]{KWitt}, that
\begin{displaymath}
\TR^{n-1}_*(\mathbb{F}_p; p) =
S_{\mathbb{Z}/p^{n-1}}\{\sigma_{n-1}\}
\end{displaymath}
where $|\sigma_{n-1}| = 2,$ and  $R(\sigma_{n-1}) =
p\lambda_{n-1}\sigma_{n-2},$ where $\lambda_{n-1} \in
\mathbb{Z}/p^{n-1}$ is a unit. It follows from \cite[Theorem B (iii)]{closure} that this unit $\lambda_{n-1} = 1$.

Let $\beta$ be a virtual representation such that $\beta'=\alpha$. Then there is a commutative diagram
$$
\xymatrix{
 \TR^n_{\alpha}(\mathbb{F}_p;p) \ar[rr]^-{\sigma_n} \ar[d]^{\hat{\Gamma}_n} & & \TR^n_{2+\alpha}(\mathbb{F}_p;p) \ar[d]^{\hat{\Gamma}_n} & \\
\hat{\mathbb{H}}^{0}(C_{p^n}, T(\mathbb{F}_p) \wedge S^{-\beta})
\ar[rr]^-{\hat{\Gamma}_n(\sigma_n)} && \hat{\mathbb{H}}^{-2}(C_{p^n},
T(\mathbb{F}_p) \wedge S^{-\beta})  &\\}
$$
Let $q= -\min(|\alpha|, |\alpha'|, |\alpha''|, \ldots |\alpha^{(n-1)}|)$. Then $\hat{\Gamma}_n$ is an isomorphism on $\TR^n_{2q + \alpha}(\mathbb{F}_p;p)$ by Proposition \ref{GammaHat}. If we understand the maps given by multiplication by $\sigma_n$ we can use the above commutative diagram to understand the map $\hat{\Gamma}_n$ by inducting down from  $\TR^n_{2q + \alpha}(\mathbb{F}_p;p)$ to $\TR^n_{\alpha}(\mathbb{F}_p;p)$.

\begin{proposition}\label{BigR}
Let $r_n = \min(n, |\alpha| + 1)$. If $r_n \geq 1$ then the maps
$$
\hat{\Gamma}_n, \Gamma_{n-1}: \TR^n_{\alpha}(\mathbb{F}_p;p) \rightarrow \mathbb{Z}/p^n
$$
agree up to multiplication by a unit.  If $r_n < 1$, then the map $\hat{\Gamma}_n: \TR^n_{\alpha}(\mathbb{F}_p;p) \rightarrow \mathbb{Z}/p^n$ is given as the composite  $\phi_n\hat{\Gamma}_{n-1}$ where $\phi_n:\mathbb{Z}/p^{n-1} \rightarrow \mathbb{Z}/p^{n}$, $\phi_n(1) \doteq p^{1-r_n}$
\begin{proof}
By the five lemma, in the range where $\hat{\Gamma}_{n-1}$ is an isomorphism, the map $\Gamma_{n-1}$ is as well. So for $q = -\min(|\alpha|, |\alpha'|, |\alpha''|, \ldots |\alpha^{(n-1)}|)$ we get a commutative diagram
$$
\xymatrix{
 \TR^n_{\alpha}(\mathbb{F}_p;p) \ar[rr]^-{(\sigma_n)^q} \ar[d]^{{\Gamma}_{n-1}} & & \TR^n_{2q+\alpha}(\mathbb{F}_p;p) \ar[d]^{{\Gamma}_{n-1}} & \\
{\mathbb{H}}^{0}(C_{p^n}; T(\mathbb{F}_p) \wedge S^{-\beta})
\ar[rr]^-{({\Gamma}_n(\sigma_n))^q} && {\mathbb{H}}^{-2q}(C_{p^n};
T(\mathbb{F}_p) \wedge S^{-\beta})  &\\}
$$
where the map $\Gamma_{n-1}$ on the right is an isomorphism. Our aim is to use this diagram to compute multiplication by $(\sigma_n)^q$. Suppose that on the left $\Gamma_{n-1}(1) \doteq p^a$. In order to compute the top horizontal map, we need to first understand the bottom horizontal map. It follows from $R(\sigma_n) = p \sigma_{n-1}$ that multiplication by $\Gamma_n(\sigma_n)$ from ${\mathbb{H}}^{-2s}(C_{p^n}; T(\mathbb{F}_p) \wedge S^{-\beta})$ to ${\mathbb{H}}^{-2s - 2}(C_{p^n}; T(\mathbb{F}_p) \wedge S^{-\beta})$ is multiplication by $p$ (up to units) for $s < 1-r_n$, and an isomorphism for $ s\geq 1 - r_n$. Consider first the case where $r_n \geq 1$. Then, $s \geq 1-r_n$ for all multiplication by $\Gamma_n(\sigma_n)$ maps in our diagram above. Hence the composite $(\Gamma_n(\sigma_n))^q$ is an isomorphism. Since the right hand $\Gamma_{n-1}$ map is also an isomorphism, this implies that $(\sigma_n)^q(1) \doteq p^a$. From the commutative diagram relating $\hat{\Gamma}_n$ and multiplication by $(\sigma_n)^q$, we then conclude that $\hat{\Gamma}_n(1)\doteq p^a$. Hence the maps $\Gamma_{n-1}$ and $\hat{\Gamma}_n$ agree up to multiplication by a unit.

We now consider the case where $r_n <1$. Since multiplication by $\Gamma_n(\sigma_n)$ is (up to units) multiplication by $p$ for $s < 1-r_n$, multiplication by $(\Gamma_n(\sigma_n))^q$ is given by  $(\Gamma_n(\sigma_n))^q(1) \doteq p^{1-r_n}.$ Since the right hand $\Gamma_{n-1}$ map is an isomorphism, this implies that $(\sigma_n)^q(1) \doteq p^{a+1-r_n}$. From the commutative diagram relating $\hat{\Gamma}_n$ and multiplication by $(\sigma_n)^q$, we then conclude that $\hat{\Gamma}_n(1)\doteq p^{a+1-r_n}$. This completes the proof.

\end{proof}
\end{proposition}

\section{Computation of $\TR^n_{\alpha}(\mathbb{F}_p;p)$} \label{Proofs}

In this section we prove Theorems \ref{main} and \ref{main2}. We start by proving Theorem \ref{main}
\begin{proof}
We prove this by induction. The result for $\TR^1_{\alpha^{(n-1)}}(\mathbb{F}_p;p)$ and the map $\hat{\Gamma}_1$ follow from our computation at the beginning of Section \ref{compute}. Now assume we have computed $\TR^{j-1}_{\alpha^{(n-j+1)}}(\mathbb{F}_p;p)$ and the map $\hat{\Gamma}_{j-1}$. We consider two cases: $r_j \geq 1$, and $r_j <1$. If $r_j \geq 1$, it follows from Proposition \ref{Pullback1} that $\TR^j_{\alpha^{(n-j)}}$ is given by the pullback
$$
\xymatrix{
\TR^j_{\alpha^{(n-j)}}(\mathbb{F}_p;p) \ar[r]^{R} \ar[d]^{{\Gamma}_{j-1}} &\TR^{j-1}_{\alpha^{(n-j+1)}}(\mathbb{F}_p;p) \ar[d]^{\hat{\Gamma}_{j-1}} \\
\mathbb{Z}/p^j \ar[r]^{R^h_j} & \mathbb{Z}/p^{j-1} \\}
$$
In order to continue inductively, we need to also compute the map $\hat{\Gamma}_j$. By Proposition \ref{BigR} the maps $\hat{\Gamma}_j$ and $\Gamma_{j-1}$ agree. Since the map $\Gamma_{j-1}$ is determined by the above pullback diagram, $\hat{\Gamma}_j$ is also determined.
We now consider the case where $r_j <1$. By Proposition \ref{Pullback1}, $\TR^j_{\alpha^{(n-j)}}(\mathbb{F}_p;p) \cong \TR^{j-1}_{\alpha^{(n-j+1)}}(\mathbb{F}_p;p)$. Proposition \ref{BigR} determines that the map $\hat{\Gamma}_j$ is as stated in the proposition.
\end{proof}
We now prove Theorem \ref{main2}
\begin{proof}
Combining the results of Proposition \ref{pullback} and Proposition \ref{BigR} , we get the following recursive definition of the length $l_{1,n} \ldots, l_{n,n}$:
Let
$$
 l_{1,1}  = \left\{\begin{array}{ll} 1  & \textup{  if    }  |\alpha^{(n-1)}| \geq 0 \\
0 &  \textup{  else  } \end{array} \right.
$$
and set $k_{1,1} =  0$.

Let $j=2$. We compute the values $l_{m,j}$ as follows, incrementing $j$ by $1$ after each iteration and stopping after the iteration where $j=n$.

Define $r_j = \min(j, |\alpha^{(n-j)}| + 1)$ and $g_{i, j-1} = l_{i, j-1} - k_{i, j-1}$  for integers $1\leq i \leq j-1$. Let $w=\#\{ 1 \leq i \leq j-1 | k_{i,j-1} \geq r_j\}$.  We choose a permutation $\tau \in \Sigma_{j-1}$ as in the statement of Proposition \ref{pullback}.  If $w = j-1$ we define $g_{\tau(w+1),j-1} = 0$.

We now give a recursive definition of the values $l_{i,j}$ and $k_{i,j}$. We set
$$
l_{1, j} = \left\{ \begin{array}{ll} \max(0, r_j) & \textup{  if  } w=j-1 \\
 \min(j, r_j+ g_{\tau(w+1), j-1}) & \textup{ else}  \\ \end{array} \right.
$$
and $k_{1,j} = 0$. For each integer $1 \leq m \leq w$, set $l_{m+1, j} = l_{\tau(m), j-1}$ and $k_{m+1, j} =  \min(l_{m+1,j}, k_{\tau(m), j-1} - r_j)$. For each integer $w+2 \leq m \leq j-1$, let $\tau(z)$ be the smallest number such that  $w+1 \leq z <m$. Then set:
$$
l_{m,j}= \left\{ \begin{array}{ll} l_{\tau(m), j-1} & \textup{  if } \tau(v) < \tau(m) \textup{ for some } w+1 \leq v <m  \\
g_{\tau(m), j-1} + k_{\tau(z), j-1} & \textup{ else} \\ \end{array}\right.
$$
Then we set $k_{m, j} = l_{m, j}$. Finally, $l_{j, j} = k_{j,j} = 0$ if $w <j-1$. If $j<n$ we increment $j$ by $1$ and iterate.

\end{proof}
The explicit algorithm given in the proof above can be easily programmed into Mathematica, providing an efficient way to compute $\TR^n_{\alpha}(\mathbb{F}_p;p)$ from the sequence of integers $|\alpha|, |\alpha'|, \ldots |\alpha^{(n-1)}|$.

We conclude by looking at a few special cases.
\subsection{Special cases}

We now interpret the main theorem in a few special cases.

\begin{corollary}
Let $\alpha$ be a virtual representation of the form $\alpha = q + \beta $, $q \in \mathbb{Z}$. Then
$$
  \TR^{n}_{\alpha}(\mathbb{F}_p;p) \cong \mathbb{Z}/p^{L}
$$
where $L$ is given as follows:
Let $m$ be such that $ |\alpha^{(n-m)}| > 0 \geq |\alpha^{(n-m+1)}|$. Let $r_j = \min(j, |\alpha^{(n-j)}| + 1)$. Set $l_{1, m} = r_m$. Letting $j$ range from $m+1$ to $n$,
$$
l_{1,j} =  \min(j, l_{1, j-1} + r_j).
$$
Then $L = l_{1, n}$.
\begin{proof}
For a representation $\alpha$ of this form,
$$
|\alpha| \geq |\alpha'| \geq \ldots \geq |\alpha^{(n-1)}|
$$
and given the hypotheses of the corollary,
$$
r_1 \leq r_{2} \leq \ldots r_{m -1} \leq 0 < r_{m} \leq \ldots r_n.
$$
From the recursive definition in the proof of Theorem \ref{main2} we know that for $j<m$, all $l_{i, j} = 0$, $1 \leq i \leq j$. It follows that $l_{1, m} = r_m$, $k_{1, m} = 0$, and $l_{i,m} = k_{i,m} = 0$ for $2 \leq i \leq m$. We show by induction that for all $j$ from $m$ to $n$, $l_{1, j} \neq 0$, $k_{1,j} = 0$, and $l_{i, j} = k_{i,j} = 0$ for $2 \leq i \leq j$. We have already noted that this holds for $j=m$. If it holds for $j$ it holds for $j+1$ by the recursive definition. Thus $\TR^n_{\alpha}(\mathbb{F}_p;p)$ is a cyclic group. The order of this group as stated in the corollary follows directly from the recursive definition.
\end{proof}
\end{corollary}

\begin{corollary}
Let $\alpha$ be a representation of the form $\alpha = q- \gamma$, $q \in \mathbb{Z}$. Then
$$
  \TR^{n}_{\alpha}(\mathbb{F}_p;p) \cong \mathbb{Z}/p^m
$$
where $m$ is the number such that
$$
|\alpha^{(n-m)}| \geq 0 > |\alpha^{(n-m - 1)}|.
$$

\begin{proof}
For a representation $\alpha$ of this form,
$$
|\alpha^{(n-1)}| \geq |\alpha^{(n-2)}| \geq \ldots  \geq |\alpha|
$$
and
$$
r_n \leq r_{n-1} \leq \ldots \leq r_{m+1} \leq 0 < r_m \leq \ldots \leq r_1.
$$

Following the recursive definition in the proof of Theorem \ref{main2},  $l_{1,1} = 1$ and $k_{1,1}=0.$ We show inductively that $l_{1,j} = j$ and $k_{1,j} = 0$ for $1 \leq j \leq m$. We have just observed that this holds for $j=1$. Now suppose it holds for $j-1$. Then by the recursive definition
$$
l_{1, j} = \min(j, r_j + g_{\tau(w+1), j-1}) = \min(j, r_j + l_{1, j-1}) = j.
$$
It is also clear that $l_{i, j} = 0$ for all $2 \leq i \leq j$. So $l_{1,m} = m$ and $l_{i,m} =0$ for $2 \leq i \leq m$. Continuing the recursion, we see that the group is unchanged for $j>m$, since in those cases $w=j-1$, and $r_j <0$ so the summand $l_{1,j} = 0$. Thus $\TR^n_{\alpha}(\mathbb{F}_p;p) \cong \mathbb{Z}/p^m$.
\end{proof}
\end{corollary}
The second corollary agrees with the result given by Hesselholt and Madsen in this case ~\cite[Proposition 9.1]{KWitt}. We have now seen that if $\alpha$ is an actual complex representation of $S^1$, or zero minus an actual representation, the $\TR$--groups we are calculating are all cyclic.

%
%
%
%

\end{document}